\documentclass{amsart}

\usepackage{amssymb}
\usepackage{amsmath}
\usepackage{amsthm}
\usepackage[mathscr]{eucal}

\newtheorem{thrm}{Theorem}[section]
\newtheorem{lemma}[thrm]{Lemma}
\newtheorem{prop}[thrm]{Proposition}

\theoremstyle{definition}

\theoremstyle{remark}

\numberwithin{equation}{section}

\newcommand{\dbar}{$\bar{\partial}$}

\begin{document}

\bibliographystyle{plain}

\title[Boundary value problems]{Boundary value problems on product domains}

\author{Dariush Ehsani}

\address{Department of Mathematics, Penn State University - Lehigh Valley, Fogelsville,
PA 18051}
 \email{ehsani@psu.edu}


\begin{abstract}
We consider the inhomogeneous Dirichlet problem on product
domains.  The main result is the asymptotic expansion of the
solution in terms of increasing smoothness up to the boundary.  In
particular, we show the exact nature of the singularities of the
solution at singularities of the boundary by constructing singular
functions which make up an asymptotic expansion of the solution.
\end{abstract}
\date{}
\maketitle
\section{Introduction}
In \cite{Eh07}, we examined the \dbar-Neumann problem for
$(0,1)$-forms on product domains in $\mathbb{C}^n$ of the form
$D=D_1\times\cdots\times D_n$ where $D_i\subset\mathbb{C}$.  In
our analysis we related the solution of the \dbar-Neumann problem
to the inhomogeneous Dirichlet problem, and as a corollary
obtained specific information on the nature of the singularities.
We were able to construct singular functions which we used to
write an asymptotic expansion of the solution, each successive
term in the expansion exhibiting a higher degree of
differentiability up to the boundary.  This paper is written to
generalize this corollary result of \cite{Eh07}.  Although we
focus on the inhomogeneous Dirichlet problem, the methods used can
be applied to other boundary value problems as well.

 We let
$\Omega\subset\mathbb{R}^n$ be a product of $q$ smooth bounded
domains, $\Omega=\Omega_1\times\cdots\times\Omega_q$, where
$\Omega_i\subset\mathbb{R}^{1+j_i}$, for $1\le j_i$ and
$\sum_{i=1}^q j_i = n-q$.  Locally, in a neighborhood, $U$ of
$x_0$, a point in the distinguished boundary of $\Omega$,
$\partial\Omega_1\times\cdots\times\partial\Omega_q$, $\rho_i$ is
a defining function for $\Omega_i$, which we shall assume to be of
the form
\[\rho_i=\phi_i(t_1^i,\ldots,t_{j_i}^i)-x_i \]
for $1\le i\le q$, and $\phi_i\in C^{\infty}(\mathbb{R}^{j_i})$.
Thus, $\Omega\cap U$ is the set of all $x\in\Omega$ such that
$\rho_i(x)<0$ for all $1\le i\le q$.

 We consider
the inhomogeneous Dirichlet problem on the product domain,
$\Omega$
\begin{align}
 \label{inhomo}
\triangle u&=f \quad \mbox{ in } \Omega \\
\nonumber
 u&=0 \quad \mbox{ on } \partial\Omega
\end{align}
and the singular behavior of $u$ in a neighborhood, $U$, of a
point, $x_0$, of the distinguished boundary, at which
$\partial\Omega$ is not smooth. In our analysis $f$ will be in the
class of $C^{\infty}(\overline{\Omega})$.  Let
$H^{\alpha}(\Omega)$ denote the Sobolev $\alpha$ space.  Existence
and uniqueness of a solution follows from Jerison-Kenig
\cite{JK81}:
\begin{thrm} [Jerison-Kenig]
 \label{j-k}
Let $\Omega$ be a bounded Lipschitz domain in $\mathbb{R}^n$.
Suppose that $\frac{1}{2} < \alpha < \frac{3}{2}$.  Then the
inhomogeneous Dirichlet problem (\ref{inhomo}) has a unique
solution $u\in H^{\alpha}(\Omega)$ and
\[
 \|u\|_{H^{\alpha}(\Omega)} \lesssim \|f\|_{H^{\alpha-2}(\Omega)}
\]
for every $f\in H^{\alpha-2}(\Omega)$.
\end{thrm}
In the same paper \cite{JK81} (see also \cite{JK95}) Jerison and
Kenig also proved that in our situation, $u\in
H^{\frac{3}{2}}(\Omega)$. In contrast to the case of smooth
boundary, the non-smooth case exhibits singularities.  Thus the
classical $L^2$ treatment (see Lions-Magenes \cite{LM}) in which
gains in derivatives are obtained has to be modified.  In this
paper it is our purpose to write an explicit solution as a sum of
terms with increasing degrees of differentiability up to the
boundary, and thus give an analysis of the behavior of the
solution near the singular parts of the boundary.

In \cite{Eh07}, the case of $\Omega_j\in\mathbb{R}^2$ for each $j$
was considered and a conformal mapping was used to reduce the
problem to one on the product of half-planes.  Instead of a
conformal mapping, in this paper we use a change of coordinates,
so that, locally, the domain looks like a product of half-spaces.

We organize our paper as follows.  In Section 2, we setup the
problem, and transform it to one on a product of half-planes.  We
begin the process of localizing the problem, so that we only
concern ourselves with function behavior in a neighborhood of a
singular point of the boundary. In Section 3 we write down an
infinite sum, which represents a solution modulo terms smooth up
to the boundary for the problem setup in Section 2. In Section 4
we construct explicit singular functions which comprise the terms
in the asymptotic expansion of our solution.

\section{A problem on a product of half-planes}
\label{change}

Recall from above that
$\Omega=\Omega_1\times\cdots\times\Omega_q$, and $\Omega_i$ in a
neighborhood of a point on the distinguished boundary has as a
defining function
 \begin{equation*}
\phi_i(t_1^i,\ldots,t_{j_i}^i)-x_i.
 \end{equation*}
   We thus use a
transformation of coordinates
\begin{align}
\label{transform}
y_j=&x_j-\phi_j \qquad 1\le j\le q\\
\nonumber y_{q+m}=&t_{k}^l
 \qquad 1\le m\le n-q, \quad j_1+\cdots+j_{l-1} +k =m, \quad
k\le j_l
\end{align}
where $x_0$ corresponds to $y=0$,
 and the related matrix $A=[a_{ij}]$, where
\begin{equation*}
a_{ij}=\begin{cases}
 1 & \quad i=j\\
 -\phi_{j,k} &\quad i\ge q, j=l\\
 0 & \quad i\ne j, j>q,
\end{cases}
\end{equation*}
and $\phi_{j,k}=\frac{\partial \phi_j}{\partial t^j_k}$, where $k$
and $l$ are determined by the unique representation of $i$ as
$i=j_1+\ldots+j_{l-1}+k$, and $k\le j_l$.

The transformation (\ref{transform}) leads to a Dirichlet problem
on the domain $\mathbb{H}^n_q$ in $\mathbb{R}^n$ given by
$\mathbb{H}^n_q=\mathbb{H}^1\times\cdots\times\mathbb{H}^q\times\mathbb{R}^{n-q}$,
where $\mathbb{H}^i=\{(y_1,\ldots,y_n):y_i>0\}$ for $i=1,\ldots,q$
in which the differential operator, $\triangle$, is replaced with
\begin{align*}
\triangle' &= \sum_{k=1}^{n}\left(\sum_{i=1}^n
a_{ki}\frac{\partial}{\partial y_i}\right)
 \left(\sum_{j=1}^na_{kj}\frac{\partial}{\partial y_j}\right)\\
 &=\sum_{i,j=1}^ng^{ij}
 \frac{\partial^2}{\partial y_i\partial
 y_j}+ \sum_{k=1}^n\sum_{i,j=1}^n
  a_{ki}\frac{\partial a_{kj}}{\partial y_i}\frac{\partial}{\partial y_j}
%
\end{align*}
where $g^{ij}=\sum_{k=1}^n a_{ki}a_{kj}$ is the metric tensor
given by the $(i,j)$ entries of the matrix $A^tA$.  Because of the
invertibility of $A$, and from the smoothness up to the boundary
of the $\rho_i$, we see that $g^{ij}\in
C^{\infty}(\overline{\mathbb{H}}^n_q)$ as are the terms
$a_{ki}\frac{\partial a_{kj}}{\partial y_i}$. Calculating $g^{ij}$
gives
\begin{align*}
\nonumber
\triangle'=&\sum_{i=1}^q
(1+|\nabla\phi_i|^2)\frac{\partial^2}{\partial
y_i^2}\\
\nonumber
& +2\sum_{i=1}^q\sum_{k=1}^{j_i}\left(-\frac{\partial\phi_i}{\partial
 t_k^i}\right) \frac{\partial}{\partial
 y_{q+j_1+\cdots+j_{i-1}+k}}\frac{\partial}{\partial y_i}
 +\sum_{i=1}^q\sum_{j,k=1}^n
  a_{kj}\frac{\partial a_{ki}}{\partial y_j}\frac{\partial}{\partial
  y_i}\\
  \nonumber
  & +\sum_{i=q+1}^n \frac{\partial^2}{\partial y_i^2}+\sum_{i=q+1}^n\sum_{j,k=1}^n
  a_{kj}\frac{\partial a_{ki}}{\partial y_j}\frac{\partial}{\partial
  y_i}.\\
%
\end{align*}
Therefore, using the change of coordinates (\ref{transform}) we
examine the problem
\begin{align*}
 \varphi \triangle' u&=\varphi f \quad \mbox{ in }
 \mathbb{H}^n_q\\
u&=0 \quad \mbox{ on } y_i=0, \ \ i=1,\ldots,q \nonumber,
\end{align*}
where $\varphi\in C^{\infty}_0 (\overline{\Omega})$ is cutoff
function such that $0\le \varphi\le 1$ and supp $\varphi \subset
U$, and $\varphi\equiv 1$ near $x_0$.

Commuting $\varphi$ with the operator $\triangle'$ gives
\begin{equation}
\label{defineh}
 \triangle' \varphi u=\varphi f +[\triangle',\varphi]u.
\end{equation}

 We begin our discussion of the singularities of the solution by
 noting the location of the singularities must be along
 intersections of boundaries of more than one $\partial \Omega_i$.
\begin{lemma}
\label{singlemm} Let $f\in C^{\infty}(\overline{\Omega})$, $u\in
L^2(\Omega)$ be the unique solution to the Dirichlet problem,
(\ref{inhomo}), on $\Omega$. Let $V\subset \Omega$ such that
$V\cap\partial\Omega_j\neq \emptyset$ for at most one $j$. Then
$u\in C^{\infty}(\overline{V})$.
\end{lemma}
\begin{proof}
The lemma is a consequence of regularity at the boundary of the
Dirichlet problem \cite{Fo}.
\end{proof}

\section{An asymptotic construction of a solution}
         We set $v=\varphi u$ and $h$ to be the right hand side of
(\ref{defineh}), and note that $h\equiv f$ in a neighborhood of
$x_0$.  Without loss of generality we can assume $h$ is just
$\varphi f$ since the error term $[\triangle',\varphi]u$ leads to
terms smooth up to the boundary in a small neighborhood of our
chosen point $x_0$. This can be seen by writing the solution in
terms of Green's function.  Let $\phi(y)$ be a function which is
equivalently 0 near $x_0$, and let $G(x,y)$ be the Green's
function for the Dirichlet problem on $\Omega$, then
$\int_{\Omega}G(x,y) \phi(y) dV(y)$ can be extended smoothly
across the boundary near $x_0$.

 We shall make use of odd reflections along $y_i=0$ for
$1\le i\le q$.  We take odd reflections of (\ref{defineh}) with
respect to the variables $y_i$ for $1\le i\le q$ and we denote the
extension by a superscript $o$.  We also denote by a superscript
$o_j$ an even extension in variables $y_i$ for $1\le i\le q$,
$i\ne j$ and an odd extension in $y_j$, and a superscript $e$ is
used to denote even extensions in all variables $y_i$ for $1\le
i\le q$.  We write (\ref{defineh}), after reflections, in the form
\begin{equation}
 \label{lapodd}
\sum_{i=1}^q \left(a_i^e \frac{\partial^2}{\partial
y_i^2}+b_i^{o_i}\frac{\partial}{\partial y_i}+c_i^{e}\right)
 v^o=h^o.
\end{equation}
 Due to our discussion above on the error term in
(\ref{defineh}) leading to terms which are $C^{\infty}$ smooth up
to the boundary of $\Omega$ in a neighborhood of $x_0$, we look to
solve (\ref{lapodd}) with $h^o$ replaced by $(\varphi f)^o$.

 We make
the note here that the difficulty in following the process in the
case of smooth domains to attempt to construct a parametrix lies
in the existence of singularities in the symbol along $y_i=0$ for
$1\le i\le q$ due to the reflections.   This difficulty can be
resolved by referring to Lemma \ref{singlemm}, and we could then
proceed to show the error terms of the parametrix construction
yield an error term in the solution which is also smooth up to the
boundary as we do in Theorem \ref{asympthrm} below.  Since this
examination of the error terms is independent of the method of
parametrix construction, we choose here to follow our analysis in
\cite{Eh07} from the beginning, which nonetheless has the flavor
of a parametrix construction.

  Here  $\eta_i$ will be the
Fourier variable corresponding to $y_i$ for $i\le q$, and $\xi_i$
the Fourier variable corresponding to $y_{i}$ for $q+1\le i\le n$.
We denote the symbols of $b_i$ and $c_i$ by $B_i$ and $C_i$,
respectively, for $1\le i\le q$.  Let $(q,i)=q+j_1+\cdots+j_{i-1}$
and from the symbol of the operator in (\ref{lapodd}),
\begin{equation*}
P(y,\xi,\eta)=-\sum_{i=1}^q (1+|\nabla\phi_i|^2) \eta_i^2+
    \sum_{i=1}^qB_i^{o_i}(y_{q+1},\ldots,y_n,\xi)
 \eta_i-\sum_{i=1}^qC_i(y_{q+1},\ldots,y_n,\xi),
\end{equation*}
we define the operator $K$ by
\begin{equation*}
Ku=- \sum_{i=1}^q\int \left(B_i^{o_i}\eta_i-C_i\right)
 \hat{u} e^{iy\cdot(\xi,\eta)} d\xi d\eta,
\end{equation*}
and with
\begin{equation}
\label{v0}
v_0=\frac{1}{(2\pi)^q}\int\chi(\eta)\frac{\widetilde{(\varphi
f)}^o
  (\eta,y_{q+1},\ldots,y_n)}
{-\sum_{i=1}^q (1+|\nabla\phi_i|^2)\eta_i^2}
 e^{i(y_1,\ldots,y_q)\cdot\eta} d\eta,
\end{equation}
where \ $\widetilde{}$ \  refers to a partial Fourier transform in
the variables $y_1,\ldots,y_q$,
 we inductively define
\begin{equation}
 \label{defnj}
v_{j+1}=\frac{1}{(2\pi)^q}\int\chi(\eta)\frac{\widetilde{Kv_j}
  (\eta,y_{q+1},\ldots,y_n)}
{-\sum_{i=1}^q (1+|\nabla\phi_i|^2)\eta_i^2}
 e^{i(y_1,\ldots,y_q)\cdot\eta} d\eta \qquad j\ge 0.
\end{equation}

We write the solution $v^o$ to (\ref{lapodd}) in the form
\begin{equation}
\label{vexpn}
 v^o= v_0+\cdots+ v_N +v_{R_N}
\end{equation}
for any integer $N\ge 0$ where the remainder term, $v_{R_N}$,
satisfies
\begin{multline}
 \label{vrn}
\sum_{i=1}^q \left(a_i^e \frac{\partial^2}{\partial
y_i^2}+b_i^{o_i}\frac{\partial}{\partial y_i}+c_i^{e}\right)
 v_{R_N} = \\Kv_N
+ F.T.^{-1}_{\eta} \left((1-\chi) \left(\sum_{j=1}^{N-1}
\widetilde{Kv_j}
 +\widetilde{(\varphi f)}^o\right)\right).
\end{multline}
We write $F.T.^{-1}_{\eta}$ to refer to an inverse Fourier
transform in the $\eta$ variables only.

   The forms of the $v_j$ are written in the following lemma.
\begin{lemma}
 \label{paraj} Each $v_j$ is of the form
\begin{equation}
 \label{formv}
v_j=
 \int_{\mathbb{R}^{n}}
 \chi^{j+1}(\eta) \frac{p_j(y_{q+1},\dots,y_n,\xi,\eta)}
 {(\sum_{i=1}^q (1+|\nabla\phi_i|^2)\eta_i^2)^{3j+1}}
 \widehat{(\varphi f)}^o e^{iy\cdot (\xi,\eta)} d\xi d\eta
\end{equation}
where  $p_j(y,\xi,\eta)$ is $C^{\infty}$ smooth as a function of
$y_{q+1},\dots,y_n$ and is a polynomial in $\xi$ and $\eta$, of
order $5j$ in the $\eta$ variables.
\end{lemma}
\begin{proof}
The lemma would be clear if it were not for the odd reflections in
the symbols, $B_k$.  We thus examine terms of the form $
H(y_k)\eta_k \tilde{u}^o,$ where $H(t)$ is defined to be 1 for
$t\ge0$ and $-1$ for $t<0$.
 We have
\begin{equation*}
F.T.^{-1}_{\eta} \left(\eta_k\tilde{u}^o\right) =i\frac{\partial
u}{\partial y_k}^{e_k} + i\delta(y_k)u^o,
\end{equation*}
 where $e_k$ denotes an even reflection in the $y_k$ variable, and
 odd reflections in\\ $y_1,\ldots, y_{k-1},y_{k+1},\ldots,y_q$.
 And thus
\begin{equation*}
H(y_k) F.T.^{-1}_{\eta} \left(\eta_k\tilde{u}^o\right)
=i\left(\frac{\partial u}{\partial y_k}\right)^{o} +
i\delta(y_k)u^o.
\end{equation*}
Using
\begin{equation*}
\widetilde{\left(\frac{\partial u}{\partial y_k}\right)^{o}}
 =-i\eta_k \widehat{u}^{e_k}
\end{equation*}
 in (\ref{defnj}), and noting that, since $\hat{v}_j$ is odd in
 all $\eta$ variables by our construction, $v_j(y_k=0)=0$ for
 $k=1,\ldots, q$,
we can write
\begin{align}
 \nonumber
-\sum_{k=1}^q (1+|\nabla\phi_k|^2)\eta_k^2 \tilde{v}_{j+1}  &=
\sum_{k=1}^q \Big( B_k(y_{q+1},\ldots, y_n,\xi)
 \left(i\eta_k \tilde{v}_j ^{e_k}+ iF.T._{y_{\bar{k}}}
 v_j(y_k=0)\right)\\
 \nonumber
&\qquad \qquad + C_k(y_{q+1},\ldots, y_n,\xi) \tilde{v}_j\Big)\\
\label{inducte}
 &= \sum_{k=1}^q \Big( B_k(y_{q+1},\ldots, y_n,\xi)
 i\eta_k \tilde{v}_j ^{e_k}
+ C_k(y_{q+1},\ldots, y_n,\xi) \tilde{v}_j\Big),
\end{align}
where $F.T._{y_{\bar{k}}}$ refers to a Fourier transform in all
$y_1,\ldots,y_q$ variables but $y_k$.  Now
\begin{equation*}
-\sum_{k=1}^q (1+|\nabla\phi_k|^2)\eta_k^2 (-i\eta_l \tilde{v}_j
^{e_l})=-i\eta_l \chi(\eta)\widetilde{Kv_{j-1}}^{e_l} +
 i\eta_l F.T._{y_{\bar{k}}} \frac{\partial v_j}{\partial
 y_l}(y_l=0)
 ,
\end{equation*}
and so
\begin{equation}
 \label{ej}
\tilde{v}_j ^{e_l}=-
\frac{\chi(\eta)\widetilde{Kv_{j-1}}^{e_l}}{\sum_{k=1}^q
(1+|\nabla\phi_k|^2)\eta_k^2} +
 \frac{F.T._{y_{\bar{k}}} \frac{\partial v_j}{\partial
 y_l}(y_l=0)}{\sum_{k=1}^q
(1+|\nabla\phi_k|^2)\eta_k^2}
  \qquad j\ge1,
\end{equation}
while for $v_0$ we have
\begin{equation}
\label{e0} \tilde{v}_0 ^{e_l}=-
\frac{\chi(\eta)\widetilde{(\varphi f)}^{e_l}}{\sum_{k=1}^q
(1+|\nabla\phi_k|^2)\eta_k^2} +
 \frac{F.T._{y_{\bar{k}}} \frac{\partial v_0}{\partial
 y_l}(y_l=0)}{\sum_{k=1}^q
(1+|\nabla\phi_k|^2)\eta_k^2}
\end{equation}

$v_0$, from (\ref{v0}), has the form of (\ref{formv}) and
  taking an inverse transform
 with respect to $\eta_l$ of the last
term on the right side of (\ref{e0}) shows that it also has a
Fourier transform in $y_1,\ldots,y_q$ of the form in
(\ref{formv}).  We use an induction argument and assume
$\tilde{v}_j$ is of the form in (\ref{formv}).

  The proof will be complete, by (\ref{inducte}), if
we show $v_j^{e_k}$ is also of the form (\ref{formv}).  An inverse
transform with respect to $\eta_k$ of the last term on the right
side of (\ref{ej}) shows that it also of the form in
(\ref{formv}).
 For the first term on
the right of (\ref{ej}) we use
\begin{align*}
 \widetilde{Kv_{j-1}}^{e_l}
=& - \sum_{i=1\atop i\ne l}^q F.T._{\eta} \int
\left(B_i^{o_i}\eta_i-C_i\right)
 \hat{v}_{j-1}^{e_l} e^{iy\cdot(\eta,\xi)} d\eta d\xi\\
&- \left(B_l\eta_l-C_l\right)\hat{v}_{j-1}\\
=& - \sum_{i=1\atop i\ne l}^q F.T._{\eta} \int
\left(B_i\eta_i-C_i\right)
 \hat{v}_{j-1}^{e_{li}} e^{iy\cdot(\eta,\xi)} d\eta d\xi\\
&- \left(B_l\eta_l-C_l\right)\hat{v}_{j-1} ,
\end{align*}
 where $e_{li}$ means an even extension in both $y_l$ and $y_i$,
 and we use the fact that an even extension in $y_l$ of a
function $H(y_l) F.T.^{-1}_{\eta} \left(\eta_l\tilde{u}^o\right)$
restricted to $y_l>0$ has a transform $\eta_l\tilde{u}^o$ if
$u(y_l=0)=0$.  The same arguments above showing $v_0^{e_k}$ is of
the form (\ref{formv}) shows we may apply the induction hypothesis
to ${v}_{j-1}^{e_{li}}$ as well.  The induction hypothesis then
shows that $\tilde{v}_j^{e_l}$ also has the form in (\ref{formv}),
and with
 $\tilde{v}_j^{e_l}$ inserted into (\ref{inducte}) we have
 $\tilde{v}_{j+1}$ of the form in (\ref{formv}).
\end{proof}

We show that (\ref{vexpn}) gives an asymptotic expansion of our
solution $v$.
\begin{thrm}
\label{asympthrm} Each successive term in (\ref{vexpn}) is of
increasing differentiability up to the boundary, and the remainder
terms $v_{R_N}$ are also of increasing differentiability up to the
boundary.  The representation (\ref{vexpn}) is therefore an
asymptotic expansion of the solution to the Dirichlet problem.
\end{thrm}
\begin{proof}
  From our definition of the $v_j$, an induction
argument shows each $v_j$ is infinitely differentiable with
respect to the variables $y_{q+1},\ldots,y_n$.  Furthermore, it
follows by our construction of the $v_j$, and by (\ref{formv}),
that if $v_j\in H^k(\mathbb{R}^n)$, the Sobolev-$k$ space, then
$v_{j+1}\in H^{k+1}(\mathbb{R}^n)$.  The proof of Lemma
\ref{paraj} also shows that restricting $v_j$ to $y_i>0$ for $1\le
i\le q$ and then reflecting about $y_i=0$ in an even manner to
form $v_j^{e_i}$ gives $v_j^{e_i}\in H^k(\mathbb{R}^n)$.  To see
the remainder terms, $v_{R_N}$ are of increasing smoothness up to
the boundary with increasing $N$, we take Fourier transforms of
(\ref{vrn}) with respect to the variables $y_1,\ldots, y_q$ to
obtain
\begin{equation}
 \label{rmexpn}
F.T._{\xi}^{-1}\left( P(y_{q+1},\ldots,y_n,\xi,\eta)
\hat{v}_{R_N}\right)
 = \widetilde{Kv_N}
+ (1-\chi) \left(\sum_{j=1}^{N-1} \widetilde{Kv_j}
 +\widetilde{(\varphi f)}^o\right).
\end{equation}
From the proof of Lemma \ref{paraj}, we have
\begin{equation*}
       \widetilde{Kv_N}
    = \sum_{k=1}^q \Big( B_k(y_{q+1},\ldots, y_n,\xi)
 i\eta_k \tilde{v}_N ^{e_k}
+ C_k(y_{q+1},\ldots, y_n,\xi) \tilde{v}_N\Big).
\end{equation*}
  If we assume
$v_N\in H^{k(N)}(\mathbb{R}^n)$, and hence $v_N^{e_i}\in
H^{k(N)}(\mathbb{R}^n)$ from above, for some $k(N)$ strictly
increasing in $N$, with the property that
$\frac{\partial^{\alpha}}{\partial y^{\alpha}} v_N\in
H^{k(N)}(\mathbb{R}^n)$ for any multi-index
$\alpha=(\alpha_1,\ldots,\alpha_n$ with
$\alpha_1=\cdots=\alpha_q=0$, then $Kv_N\in
H^{k(N)-1}(\mathbb{R}^n)$ by definition of the operator $K$, and
$p_{k(N)-1}(\eta)\widetilde{Kv_N}\in L^2$ for $p_{k(N)-1}(\eta)$ a
polynomial in the $\eta $ variables of degree $k(N)-1$.  Also,
since any polynomial of any degree in the $\eta$ variables
multiplied by the last term on the right hand side of
(\ref{rmexpn}) is in $L^2$ due to the compact support in $\eta$ of
the term $1-\chi$, we see, by multiplying (\ref{rmexpn}) by
$p_{k(N)-1}(\eta)$, that
 \begin{equation*}
p_{k(N)-1}\left(\frac{\partial}{\partial y_1},\ldots,
 \frac{\partial}{\partial y_q}\right) v_{R_N}
\end{equation*}
is the solution to a problem with mixed Dirichlet and Neumann
conditions, and by \cite{JK81}, its restriction to
$\mathbb{H}^n_q$ is in $H^{\alpha}(\mathbb{H}^n_q)$ for
$\alpha<3/2$.  Therefore, $v_{R_N}|_{\mathbb{H}^n_q}\in
H^{k(N)+\alpha}(\mathbb{H}^n_q)$ for $\alpha<1/2$, and so is of
increasing smoothness up to the boundary for increasing $N$.
\end{proof}

\section{An explicit calculation of the singularities}

In order to determine the singularities in each term in the
expansion (\ref{vexpn}), we take inverse Fourier transforms of the
expression (\ref{formv}) with respect to the $\eta$ variables.

Without loss of generality we suppose the cutoff function $\chi$
is of the form
\begin{equation*}
\chi_{\eta_1}\cdots\chi_{\eta_q},
\end{equation*}
where $\chi_{\eta_i}$ is a cutoff function in the variable
$\eta_i$ only, equal to 0 in a neighborhood of $\eta_i=0$.  We can
do this because
\begin{equation*}
\chi^{j+1}(\eta)- \chi_{\eta_1}\cdots\chi_{\eta_q}
\end{equation*}
can be written as a sum of terms, each of which has support
contained in large $\eta_k$ for at most one $k$.  Hence,
\begin{equation*}
(\chi^{j+1}(\eta)- \chi_{\eta_1}\cdots\chi_{\eta_q})
\frac{p_j(y_{q+1},\dots,y_n,\xi,\eta)}
 {(\sum_{i=1}^q a_i\eta_i^2)^{3j+1}}
 \widehat{(\varphi f)}^o,
\end{equation*}
where $a_i=1+|\nabla\phi_i|^2$, is a sum of terms which have
infinite decay in all but one $\eta_j$, and thus is the transform
of a function which is $C^{\infty}$ in all variables but one
$y_j$.  As such a term is odd in $y_j$ and as the denominator
$(\sum_{i=1}^q a_i\eta_i^2)^{3j+1}$ is the symbol of an elliptic
operator, such a term is the solution to a Dirichlet problem on a
half-plane, and is therefore $C^{\infty}$ smooth up to the
boundary of the half-plane.
 In summary the difference between (\ref{formv}) and
\begin{equation}
 \label{diffcut}
\chi_{\eta_1}\cdots\chi_{\eta_q}
\frac{p_j(y_{q+1},\dots,y_n,\xi,\eta)}
 {(\sum_{i=1}^q a_i\eta_i^2)^{3j+1}}
 \widehat{(\varphi f)}^o
\end{equation}
is the transform of a term which, when restricted to
$\mathbb{H}^n_q$, is in $C^{\infty}(\overline{\mathbb{H}}^n_q)$.

We then take (\ref{diffcut}) and integrate by parts with respect
to the variables $y_i$ for $1\le i\le q$ in the Fourier integral
of $\widehat{\varphi f}^o(\eta,\xi)$, starting with $y_1$ and
proceeding to $y_q$, keeping as remainder terms those which have
decay in one $\eta$ variable
 to the order $-2(N+q)-1$.  Such terms are of the form
 \begin{multline}
 \label{partsremainder0}
\chi_{\eta_1}\cdots\chi_{\eta_q}
\frac{p_j(y_{q+1},\dots,y_n,\xi,\eta)}
 {(\sum_{i=1}^q a_i\eta_i^2)^{3j+1}}
 \frac{1}{(i\eta)^{\alpha}}\frac{1}{(i\eta_{k+1})^{2(N+q)+1}}\times
 \\ F.T._{k+1,\ldots,n}\left( \frac{\partial^{2(N+q)+1}}{\partial
 y_{k+1}^{2(N+q)+1}} \frac{\partial^{\alpha}(\varphi f)}{\partial
 y^{\alpha}}^o (0,\ldots,0,y_{k+1},\ldots,y_n)\right)
\end{multline}
in which $\alpha$ is a $q$ index for which
$\alpha_{k+1}=\cdots=\alpha_n=0$.  We shall show below in Theorem
\ref{Dirthrm} that such remainder terms are sufficiently
continuous up to the boundary of $\mathbb{H}^n_q$.

Setting aside the remainder terms (\ref{partsremainder0}), we
analyze those other terms which result from the  expansion of
 $\widehat{(\varphi f)}^o(\eta,\xi)$, and we are led to study terms of the form
 \begin{equation}
  \label{partsremainder}
\chi_{\eta_1}\cdots\chi_{\eta_q}\frac{1}{\eta^k}\frac{\phi_{k,j}(y_{q+1},\dots,y_n,\xi,\eta)}
  {\left(\sum_{i=1}^q a_i\eta_i^2\right)^{3j+1}},
 \end{equation}
where   $\phi_{k,j}(y,\xi,\eta)$ takes the form of odd reflections
along $y_i=0$ for each $i=1\ldots q$ of functions of $y$ smooth up
to the boundary, and is a polynomial in $\xi$ and $\eta$ of order
$5j$ in the $\xi$ and $\eta$ variables, and a polynomial of degree
$j$ in the $\eta$ variables..  We use the notation
$\eta^k=\eta_1^{k_1}\cdots\eta_q^{k_q}$ for $k=(k_1,\ldots,k_q)$.

Up to multiplication by a constant the following relation holds
for $0<l<\frac{q}{2}$
\begin{equation*}
\int_{\mathbb{R}^q}\frac{1}{\left(\sum_{i=1}^q
a_i\eta_i^2\right)^l}e^{i(y_1,\cdots,y_q)\cdot\eta}d\eta =
\frac{1}{\prod_{i=1}^q \sqrt{a_i}}
 \frac{1}{\left(\sum_{i=1}^q \frac{y_i^2}{a_i} \right)^{\frac{q}{2}-l}}.
\end{equation*}

 Let
\begin{equation}
\label{Phil} \Phi_l^q=\frac{1}{\left(\sum_{i=1}^q
\frac{y_i^2}{a_i} \right)^{\frac{q}{2}-l}} \quad 0<l<\frac{q}{2}.
\end{equation}

For $l\ge \frac{q}{2}$ and $q$ even, we define $\Phi_l^q$ to be
the unique solution of the form
\begin{equation*}
p_1(y)\log\left(\sum_{i=1}^q \frac{y_i^2}{a_i} \right)+p_2(y),
\end{equation*}
where $p_2(y_1=0)=0$, $p_1$ and $p_2$ are polynomials of degree
$2l-q$ in the variables $y_i$ for $1\le i\le q$, and are
$C^{\infty}$ smooth with respect to variables $y_i$ for $q+1\le i
\le n$, to the equation
\begin{equation}
 \label{phidef}
\frac{\partial\Phi_{l}^q}{\partial y_1}=y_1\Phi_{l-1}^q.
\end{equation}
For the case $q=2$, we take
\begin{equation*}
 \Phi_{1}^2=-\frac{i}{2}
\log\left(\frac{y_1^2}{a_1}+\frac{y_2^2}{a_2}\right).
\end{equation*}

  For $q$ odd we define $\Phi_l^q$ for $l\ge
\frac{q}{2}$ as in (\ref{Phil}), the Fourier transform being
calculated in the sense of distributions.

In the sense of distributions, we take transforms of
(\ref{phidef}) and write, for $q>2$
\begin{align*}
i\eta_1\widehat{\Phi_{l}^q}&=-i\frac{\partial}{\partial \eta_1}
 \widehat{\Phi_{l-1}^q}\\
&=-2i(l-1)a_1\frac{\eta_1}{\left(\sum_{i=1}^q
a_i\eta_i^2\right)^l},
\end{align*}
which implies
\begin{equation*}
\chi_{\eta_1}\widehat{\Phi_{l}^q}=
 -2(l-1)a_1\chi_{\eta_1}\frac{1}
 {\left(\sum_{i=1}^q a_i\eta_i^2\right)^l}
\end{equation*}
for $q>2$.  For the case $q=2$, we use \cite{Eh03b} to write
\begin{equation*}
\chi_{\eta_1}\chi_{\eta_2}
 \widehat{\chi\Phi_1^2}=
 \chi_{\eta_1}\chi_{\eta_2}\frac{\sqrt{a_1a_2}}{a_1\eta_1^2+a_2\eta_2^2} +s,
\end{equation*}
where $\chi$ is a cutoff function such that $\chi\equiv1$ in a
neighborhood of the origin as in Section \ref{change}, and here
and below we use $s$ to denote terms which after
 taking inverse
 transforms give
 $C^{\infty}$ functions in some neighborhood of the origin in
$\mathbb{R}^n$.  We have therefore established the property for
$l>0$
\begin{equation*}
\chi_{\eta}
 \widehat{\chi\Phi_l^q}=
 \chi_{\eta}\frac{\prod_{i=1}^q \sqrt{a_i}}
 {\left(\sum_{i=1}^q a_i\eta_i^2\right)^l} +s,
\end{equation*}
up to multiplication by a constant.

 For $k=(k_1,\ldots,k_q)$ we further define
\begin{equation*}
\Phi_{lk}^q=\int_0^{y_q}\cdots\int_0^{t_2^{k_q}}
\cdots\int_0^{y_1}\cdots\int_0^{t_2^{k_1}}
\Phi_l^q(t_1^{k_1},\ldots,t_1^{k_q})
 dt_1^{k_1}\cdots dt_{k_1}^{k_1}\cdots dt_1^{k_q}\cdots
 dt_{k_q}^{k_q}.
\end{equation*}

We have the following key property of the $\Phi_l^q$ which allow
us to match these functions with the singularities of our infinite
sum solution to the Dirichlet problem in (\ref{vexpn}).
\begin{prop}
Let $\Phi_{lk}^q$ be defined as above.  Let $\chi\in
C^{\infty}_0(\mathbb{R}^n$ be a smooth cutoff function such that
$\chi\equiv 1$ in a neighborhood of 0.  Up to multiplication of a
smooth function of $y_i$ for $q+1\le i\le n$, $\Phi_{lk}^q$ has
the property
\begin{equation}
 \label{FTphi}
\chi_{\eta_1}\cdots\chi_{\eta_q}\widehat{\chi\Phi_{lk}^q}=
\chi_{\eta_1}\cdots\chi_{\eta_q}
 \frac{1}{\eta^k}\frac{1}{\left(\sum_{i=1}^q
a_i\eta_i^2\right)^l}+s.
\end{equation}
\end{prop}
\begin{proof}
We have
\begin{align*}
\int_{\mathbb{R}^q} \chi \Phi_{lk}^q e^{-i(y_1,\ldots,y_q)} dy
 &= \frac{1}{i^{|k|}\eta^k}\int_{\mathbb{R}^q}
 \chi \Phi_{lk}^q \frac{\partial^{|k|}}{\partial y^k}
 e^{-i(y_1,\ldots,y_q)} dy\\
&=\frac{(-1)^{|k|}}{i^{|k|}\eta^k}\int_{\mathbb{R}^q}
 \chi \frac{\partial^{|k|}}{\partial y^k} \Phi_{lk}^q
 e^{-i(y_1,\ldots,y_q)} dy +s,
\end{align*}
where the term $s$ in the last line comes from taking derivatives
of the cutoff $\chi$.  By definition
\begin{equation*}
\frac{\partial^{|k|}}{\partial y^k} \Phi_{lk}^q
 =\Phi_{l}^q.
\end{equation*}
Therefore,
\begin{equation}
 \label{lk}
\widehat{\chi\Phi_{lk}^q}=\frac{(-1)^{|k|}}{i^{|k|}\eta^k}
 \widehat{\chi\Phi_{l}^q}+s,
\end{equation}
and from the definitions and discussion above,
\begin{align}
 \nonumber
\widehat{\chi\Phi_{l}^q}&=\widehat{\Phi_{l}^q}
 +\widehat{(1-\chi)\Phi_{l}^q}\\
\nonumber
&=\widehat{\Phi_{l}^q}+s\\
 \label{justl}
&=\varphi(y_{q+1},\ldots,y_n) \frac{1}{\left(\sum_{i=1}^q
a_i\eta_i^2\right)^l}
 +s
\end{align}
for some smooth function $\varphi$.

Inserting (\ref{justl}) into (\ref{lk}) finishes the proof.
\end{proof}

For a multi-index $p=(p_1,\ldots,p_i)$, in which $p_j\le q$, we
define $\Phi_{lk}^p$ for $2\le i\le q-1$ in the same fashion as we
did $\Phi_{lk}^q$ but with respect to the $i$ variables
$y_{p_1},\ldots,y_{p_i}$, in particular $k$ is a multi-index of
length $i$.
 Thus
\begin{equation}
\label{FTphip}
 \chi_{\eta_{p}}\widehat{\chi\Phi_{lk}^p}=
\chi_{\eta_p}
 \frac{1}{\eta^k}\frac{1}{\left(\sum_{j=1}^i
a_{p_j}\eta_{p_j}^2\right)^l}+s,
\end{equation}
where $\chi_{\eta_p} = \chi_{\eta_{p_1}}\cdots\chi_{\eta_{p_i}}$.
 We also use the notation $\Phi_{lk}^p=\Phi_{lk}^q$
when $p=(1,\ldots,q)$.

From (\ref{FTphi}), (\ref{FTphip}), and the construction of the
$\Phi_{lk}^p$ we see that
\begin{equation*}
\left(\sum_{i=1}^q
a_i\eta_i^2\right)^m\eta_1^{m_1}\cdots\eta_q^{m_q}
 \left(\chi_{\eta_p}(\widehat{\chi\Phi_{lk}^p})-\widehat{\chi\Phi_{lk}^p}\right)
 \in s
\end{equation*}
for large enough $m, m_1,\ldots,m_n$.  Now
 $\chi_{\eta_p}-1$ has support for large $\eta_k$ for at most one
 $k$, and thus
$\eta_1^{m_1}\cdots\eta_q^{m_q}
 \left(\chi_{\eta_p}(\widehat{\chi\Phi_{lk}^p})-
\widehat{\chi\Phi_{lk}^p}\right)$ is the transform of a solution
to a Dirichlet problem on a half-plane with data smooth up to the
boundary.  Then, by the ellipticity of the operator whose symbol
is $\left(\sum_{i=1}^q a_i\eta_i^2\right)^m$, we obtain
derivatives of the inverse transform of
\begin{equation*}
\chi_{\eta_p}(\widehat{\chi\Phi_{lk}^p})-(\widehat{\chi\Phi_{lk}^p})
\end{equation*}
 is a $C^{\infty}$ function on a half-plane.
 By
inverting derivatives, by integration, and invoking Lemma
\ref{singlemm} to show the constants of integrations are smooth up
to the boundary, we see the terms described by
$\chi_{\eta_p}(\widehat{\chi\Phi_{lk}^p})$ and
$(\widehat{\chi\Phi_{lk}^p})$ differ by functions smooth up to the
boundary, and $\chi\Phi_{lk}^p$ will thus be seen to describe the
singularities of the solution.

With a slight abuse of notation we also use the notation
$\Phi_{lk}^p$ even after a change of variables back to the domain
$\Omega$.   We are now ready to prove the
\begin{thrm}
 \label{Dirthrm}
Let $f\in C^{\infty}(\overline{\Omega})$ and let $u\in
L^2(\Omega)$ be the unique solution to the inhomogeneous Dirichlet
problem on $\Omega$.  Then near the distinguished boundary,
$\partial\Omega_1\times \cdots\times\partial\Omega_q$, $u$ is of
the form
\begin{equation}
 \label{thrmexpn}
u=\sum_{|k|\ge0,l\ge1,\atop2\le\ell(p)\le q} c_{klp}\Phi_{lk}^p
\end{equation}
where $c_{klp}\in C^{\infty}(\overline{\Omega})$, and where
$\Phi_{lk}^p$ are defined as above.
\end{thrm}
\begin{proof}
(\ref{thrmexpn}) is obtained by following the procedure outlined
above, matching terms, (\ref{partsremainder}), in Fourier space to
the appropriate function $\Phi_{kl}^p$.  We thus need to study the
effect $\phi_{k,j}$ has as a polynomial in the $\xi$ and $\eta$
variables in (\ref{partsremainder}) on the functions
$\Phi_{kl}^p$.  We also have to treat the remainder terms
(\ref{partsremainder0}).

We first show the remainder terms given by (\ref{partsremainder0})
are described by the functions $\Phi_{lk}^p$.  We consider the
case in which (\ref{partsremainder0}) is given by
\begin{equation*}
\chi_{\eta_1}\cdots\chi_{\eta_q}
 \frac{1}{\eta_1\cdots\eta_{q-1}}\frac{1}{\eta_q^{2(N+q)+1}}
 \frac{1}{\sum_{i=1}^q
a_i\eta_i^2}F.T._{m}\left(\frac{\partial^{2(N+q)+1}}{\partial
y_m^{2(N+q)+1}}(\varphi f)^o\right)
\end{equation*}
restricted to $y_i=0$ for $1\le i\le q$ and $i\ne m$, where
$F.T._m$ denotes the partial Fourier transform with respect to
$y_m$.  The other such remainder terms are handled in a similar
manner.  The fraction
\begin{equation*}
\frac{1}{\sum_{i=1}^q a_i\eta_i^2}
\end{equation*}
is expanded in a geometric series in
\begin{equation*}
\frac{a_m\eta_m^2}{\sum_{i=1\atop i\ne m}^q a_i\eta_i^2}
\end{equation*}
up to $(N+q)/2$ terms (we assume $(N+q)$ is even), the first terms
leading to the functions $\Phi_{lk}^p$ in which $m\notin p$, while
the last term is, up to multiplication by a $C^{\infty}$ smooth
function of the variables $y_i$ for $q+1\le i\le n$, given by
\begin{equation}
 \label{solntodir}
\left(\chi_{\eta_1}\cdots\chi_{\eta_q}
 \frac{F.T._{m}\left(\frac{\partial^{2(N+q)+1}}{\partial
y_m^{2(N+q)+1}}(\varphi
f)\right)^{e_m}}{\eta_m^{N+q+1}}\frac{1}{\left(\sum_{i=1\atop i\ne
m}^q a_i\eta_i^2\right)^{(N+q)/2+1}}\right)\frac{1}{\sum_{i=1}^q
a_i\eta_i^2}.
\end{equation}
The term in parentheses is the transform of a function odd in the
variable $y_m$ and in $C^N(\mathbb{R}^{q-1}\times \{y_m>0\})$, and
so all of (\ref{solntodir}) may be viewed as the solution to a
Dirichlet problem depending $C^{\infty}$ smoothly on the
parameters $y_i$ for $q+1\le i\le n$ on a half-space, and by
regularity of the Dirichlet problem, the term in
(\ref{solntodir}), after an inverse transform, is in
$C^N(\mathbb{R}^n)$.

We finish the proof of Theorem \ref{Dirthrm} by showing that the
polynomials of the $\eta$ and $\xi$ variables in the numerator of
(\ref{partsremainder}) still preserve the form of the functions
$\Phi_{lk}^p$.  Using the fact that the data function
$(\widehat{\varphi f})^o$ has infinite decay with respect to the
$\xi$ variables, and that multiplication by $\eta_i$ corresponds
to differentiating with respect to $y_i$ the following relations
may be used to complete the proof of the theorem.  Up to
multiplication by smooth functions of $y_i$ for $q+1\le i\le n$ we
have
\begin{align*}
\frac{\partial}{\partial y_{i}}\Phi_{lk}^q&=\Phi^q_{l(k_1,\ldots,k_i-1,\ldots,k_q)},\\
\frac{\partial}{\partial y_{i}}\Phi_l^q&=y_i\Phi_{l-1}^q,
\end{align*}
with the obvious analogies for the $\Phi_{lk}^p$.

Lastly, we restrict our expansion (\ref{vexpn}) to the product of
upper half-spaces, $\mathbb{H}_n^q$, and obtain an expression of
the solution modulo terms in
$C^{\infty}(\overline{\mathbb{H}}_n^q)$ in terms of the functions
$\Phi_{lk}^p$.  After changing coordinates back to $\Omega$, we
obtain the expansion in (\ref{thrmexpn}).
\end{proof}

Our results here are comparable to those in \cite{Eh07}, in which
$\Omega_i\subset\mathbb{R}^2$, and we note that an increase in the
dimensions of the $\Omega_i$ do not affect the form of
singularities occurring. See also \cite{Eh03b} for the specific
case of $\Omega=\Omega_1\times\Omega_2\subset\mathbb{R}^4$, where
$\Omega_i\subset\mathbb{R}^2$ for $i=1,2$, in which logarithmic
and arctangent singularities are found along the corner.

Lastly, it is important to mention that there are cases in which
the solution to the Dirichlet problem does exhibit singularities.
The sum (\ref{thrmexpn}) is not trivial; the coefficients
$c_{klp}$ are not always 0.  The example $f\equiv 1$ on $\Omega$
reveals this to be the case.

\end{document}